\documentclass[10pt]{amsart}

\usepackage{amsmath,amssymb,amscd,amsthm,graphicx}
\usepackage[toc,page]{appendix}

\usepackage{hyperref}

\numberwithin{equation}{section} 
\theoremstyle{plain}
\newtheorem{thm}[equation]{Theorem}

\newtheorem{lem}[equation]{Lemma}
\newtheorem{lemma}[equation]{Lemma}
\newtheorem{cor}[equation]{Corollary}

\theoremstyle{remark}
\newtheorem{remark}[equation]{Remark}

\theoremstyle{definition}

\newtheorem{defn}[equation]{Definition}

\newtheorem*{question*}{Question}

\newcommand{\F}{{\mathcal F}}

\renewcommand{\L}{{\mathcal L}}

\newcommand{\sphere}{\ensuremath{\mathbb{S}}}
\newcommand{\reals}{\ensuremath{\mathbb{R}}}

\def\XXint#1#2#3{{\setbox0=\hbox{$#1{#2#3}{\int}$}
     \vcenter{\hbox{$#2#3$}}\kern-.5\wd0}}

\usepackage[a4paper, total={5.5in, 8.4in}, hmarginratio=1:1]{geometry}
\usepackage{hyperref}

\begin{document}

\title[The evolution of Jordan curves on $\sphere^2$]{The evolution of Jordan curves on $\sphere^2$ by \\curve shortening flow}
\author{Joseph Lauer}
\date{January 21, 2016}
\maketitle

\begin{abstract}  In this paper we prove that if $\gamma$ is a Jordan curve on $\sphere^2$ then there is a smooth curve shortening flow defined on $(0,T)$ which converges to $\gamma$ in $\mathcal{C}^0$ as $t\to 0^+ $.  Another perspective is that the level-set flow of $\gamma$ is smooth.  This is a generalization of the results of~\cite{L}, where the planar case was studied. If a Jordan curve on $\sphere^2$ has Lebesgue measure zero then we show that the level-set flow instantly becomes a smooth closed curve.  If the Lebesgue measure is positive then for small time the level-set flow is an annulus with smooth boundary.  This second case should be interpreted as a failure of uniqueness.

As in~\cite{L} key step in the proof is establishing a length estimate for smooth curves that depends on a geometric quantity called the $r$-multiplicity.  The majority of this paper concerns the extension of this length estimate to $\sphere^2$.   
\end{abstract}


\section{Introduction}  

In the study of partial differential equations a fundamental question is when, and in what sense, a solution exists for low regularity initial data.   In this paper we study the question of existence and uniqueness of curve shortening flow when the initial data is a Jordan curve on $\sphere^2$.  This result extends the author's work in~\cite{L} where the planar case was studied. The main result of this paper is the following:

\begin{thm} \label{JordanCurve} Let $\gamma:\sphere^1\to\sphere^2$ be  a  Jordan curve.  Then there exists a maximal smooth solution of curve shortening flow
$$
\gamma:\sphere^1\times(0,T)\to\sphere^2
$$
such that
$$
\lim_{t\to 0}\gamma_t=\gamma
$$
in the space of continuous curves.

Moreover, if the Lebesgue measure of $\gamma$ is zero then
\begin{enumerate}
    \item the solution is unique up to reparametrization, and
    \item $T=\infty$ if and only if the original curve bisects the area of $\sphere^2$.
\end{enumerate}
\end{thm}

As in the planar case~\cite{L} we first show that the level-set flow of $\gamma$ is smooth.  The level-set flow is a weak notion of curve shortening flow (and more generally mean curvature flow) which evolves a compact set in a way that agrees with smooth curve shortening flow when it exists.  

For positive area curves we show in Section~\ref{unique_section} that for small positive times the level-set flow is a smooth annulus which eventually either vanishes, converges to a hemisphere or takes up all of $\sphere^2$.  See Theorem~\ref{positive}.  

Let $(\Sigma,g)$ be a 2-dimensional Riemannian manifold and $\gamma_0:\sphere^1\to\Sigma$ be a smooth immersion.  A 1-parameter family of immersions $\gamma:\sphere^1\times[0,T)\to\Sigma$ is a solution to curve shortening flow with initial data $\gamma_0$ if
$$
\frac{\partial\gamma}{\partial t}=\kappa_g\vec{n},
$$
$$
\gamma(\cdot, 0)=\gamma_0,
$$
where $\kappa_g$ is the geodesic curvature and $\kappa_g\vec{n}$ is the curvature vector.    If $\gamma:\sphere^1\times(a,b)\to\Sigma$ is a solution to curve shortening flow then we denote by $\gamma_t$ the smooth curve $\gamma(\cdot, t)$.

The short-time existence of solutions for smooth initial data was proved in the planar case by Gage and Hamilton~\cite{GH86} and for surfaces which are convex at infinity by Grayson~\cite{Gray}.   In the case of embedded initial data it was also proved in~\cite{Gray} that there are only two possibilities for the long-term behaviour of such a solution.  The first is that there exists $T<\infty$ so that the solution exists only on $[0,T)$, collapses to a point and has a `circular' singularity as $t\to T$.   The second case is that the solution exists on $[0,\infty)$ and $\kappa_g$ converges uniformly to zero.  When $\Sigma=\sphere^2$ the second case occurs if and only if the original curve bisects the area of $\sphere^2$ since the Gauss-Bonnet Theorem can be used to show that this property is preserved by the flow~\cite{Gage}.  

Previously, the most general existence result in a general surface was proved by Huisken and Ecker~\cite{EH91} who required that $\gamma$ be `uniformly locally-Lipschitz', a stronger condition than rectifiability.  Recently, Hershkovits~\cite{HO} has shown that certain Reifenberg sets have smooth level-set flow, including some fractals in $\reals^{n+1}$ with $n>1$.  Thus~\cite{HO} provides the first example of such behaviour in higher dimensions.

The majority of the work towards establishing the smoothness of the level-set flow is in proving a length estimate whose statement contains the notions of a $(C,\theta)$-spacing (Section~\ref{space_sec}) and the $r$-multiplicity $M_r(\gamma)$ (Section~\ref{mult_sec}).  Roughly speaking, a $(C,\theta)$-spacing of $\gamma$ is a large collection of open balls of radius $C$ in the complement of $\gamma$.  The $r$-multiplicity acts as a coarse intersection number.  The definition is given in the outline below.

\begin{thm}\label{length}   Let $C>0$ and $\theta\ll 1$.  Then there exists $T=T(C,\theta)>0$ such that for each $0<t<T$ there exists $r=r(C,\theta,t)$ and $\widetilde{C}=\widetilde{C}(C,\theta,t)$ such that if $\gamma:\sphere^1\to\sphere^2$ is a smooth embedded curve with a $(C,\theta)$-spacing then
$$
\mathcal{L}(\gamma_t)<\widetilde{C}M_r(\gamma).
$$
\end{thm}

In this paper we prove Theorem~\ref{length} for embedded curves only.  It is possible to produce estimates for immersed curves, as was done for planar curves in~\cite{L}, but here we prove only what is necessary for the proof of Theorem~\ref{JordanCurve}.

The utility of Theorem~\ref{length} stems from the fact that if $\gamma_n\to\gamma$ uniformly then a $(C,\theta)$-spacing for $\gamma$ will be a $(C,\theta)$-spacing for $\gamma_n$ for sufficiently large~$n$.  Since the $r$-multiplicity is also well-behaved under uniform convergence we obtain the following: 

\begin{thm}\label{approx}
Let $\gamma$ be a Jordan curve on $\sphere^2$.  Then there exist constants $T, C>0$ and a function $r:(0,T)\to\reals^+$ such that if $\gamma_n$ is a sequence of smooth embedded curves that converge uniformly to $\gamma$ then
$$
\mathcal{L}\left((\gamma_n)_t\right)<CM_{r(t)}(\gamma).
$$
for $0<t<T$ and $n$ sufficiently large.
\end{thm}

The important points are that (1) the right-hand side is independent of $n$, and (2) the estimate is valid on a definite time interval independent of the approximating sequence.  Theorem~\ref{approx} immediately rules out the possibility that the level set flow of a measure zero Jordan curve has infinite length.


\subsection{Outline of the Proof of Theorem~\ref{length}} 

We begin with the $r$-multiplicity, which is a coarse intersection profile.  

\begin{defn} [$r$-multiplicity ] Let $g$ be a great circle, $0<r<\frac{\pi}{2}$ and $\gamma$ be a Jordan curve in $\mathbb{S}^2$.  Then the $r$-multiplicity of $\gamma$ at $g$, denoted by $M_{r,g}(\gamma)$, is defined as the number of components of $\gamma\cap B_{2r}(g)$ which intersect $\overline{B_r(g)}$ non-trivially. 

Moreover, we define {\it the $r$-multiplicity of $\gamma$} by
 $$M_r(\gamma)=\sup_g\{M_{r,g}(\gamma)\}.
 $$  
\end{defn}

The $r$-multiplicity has the property that for a sequence of closed curves~$\gamma_n$ which converge uniformly to a Jordan curve, the quantities $M_r(\gamma_n)$ are uniformly bounded for each $r>0$.  See Section~\ref{mult_sec}.

The plan for proving Theorem~\ref{length} is to establish local length estimates.  Let $x\in\sphere^2$.  An upper bound for the length of $\gamma_t\cap B_r(x)$ can be obtained which is proportional to the maximum number of times that $\gamma$ intersects each leaf in two transverse foliations of $B_r(x)$ with linear leaves.  This approach is unsuitable for our applications since the number of times $\gamma_n$ intersects a particular great circle is not necessarily bounded.

Instead, we replace the linear foliation by one for which the number of intersections of each leaf with $\gamma$ is controlled by the $r$-multiplicity.  These foliations are no longer linear, but evolve to be nearly linear at a given time~$t>0$.

More precisely, let $g$ be a great circle through $x$.  For sufficiently small $t>0$ we determine an appropriate scale $r>0$ and construct a foliation $\mathcal{F}$ of $B_r(g)$ such that if $\ell$ is a leaf of $\mathcal{F}$ then
\begin{enumerate}
\item $\ell_t$ is $\mathcal{C}^1$-close to $g$, and
\item $|\ell\cap\gamma|\leq 2M_r(\gamma)$.
\end{enumerate}

\begin{figure}
\centering
\scalebox{0.6}{\includegraphics{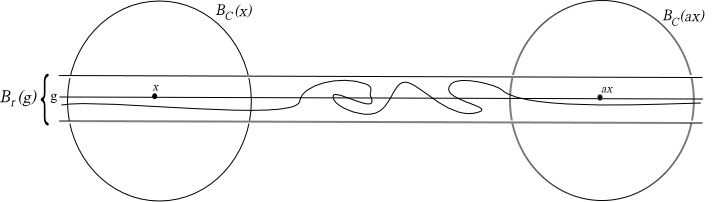}}
\caption{A leafable curve is contained in some neighborhood of a great circle $g$ and is $\mathcal{C}^1$-close to $g$ on a set of the form $B_C(x)\cup B_C(ax)$, where $a$ denotes the antipodal map on $\sphere^2$.  But outside of that set there is no restriction on the curve and in particular it need not be graphical over $g$. }  
\label{leafable_fig}
\end{figure}

In Section~\ref{straight_sec} we define the set of curves which will be allowable as leaves of~$\mathcal{F}$.  We call such curves leafable.  A leafable curve is contained in a thin neighbourhood of a great circle~$g$, but is only assumed to be $\mathcal{C}^1$-close to $g$ near two antipodal points.  In particular, leafable curves are not  necessarily graphical over $g$.  See Figure~\ref{leafable_fig}.  

The definition of leafable depends on a parameter $C>0$ which measures the size of the set on which the curve is required to be close to $g$.  In practice, $C$ is determined by the size of the spacing for the curve whose length is being bounded.  The result that guarantees that a leafable curves evolves to be nearly linear is the following:

\begin{thm}[Straightening Lemma]\label{StraightIntro} Given $C>0$ and $\alpha\ll 1$ there exists $T>0$ such that for each $0<t<T$ there exists $r=r(t,C, \alpha)>0$ with the following property:  If $g$ is a great circle and $\ell\subset B_r(g)$ is leafable then $\ell_t$ is an $\alpha$ $\mathcal{C}^1$-close to $g$.
\end{thm}

In~\cite{L} the proof of the analogous result used a family of grim reapers.  Indeed the lack of a suitable replacement is one of the main obstacles in extending the arguments in~\cite{L} to general surfaces.  For $\sphere^2$ we use a solution of the Dirichlet Problem for curve shortening flow, i.e. the evolution of an arc with fixed endpoints, to play the same role.  The existence of such an evolution is a result of Allen, Layne and Tsukhara~\cite{ALT}.

To complete the local length estimates let $\tilde{g}$ be a great circle perpendicular to $g$ at $x\in\sphere^2$ and construct a foliation as above of $B_r(\tilde{g})$.  We then show that the restriction of the two foliations at time $t$ to $B_r(x)$ is uniformly bi-Lipschitz equivalent to the standard grid in $B_r(0)\subset\reals^2$, and estimates for $\mathcal{L}(\gamma_t\cap B_r(x))$ follow since the number of intersections of $\gamma_t$ with each leaf in either foliation is at most $2M_r(\gamma)$.

In general the argument follows the planar case~\cite{L} closely.  The are two main differences: 
\begin{enumerate}
\item In~\cite{L} the initial leaves of the foliations are linear at infinity.  This is of course impossible on $\sphere^2$ so we introduce the concept of a $(C,\theta)$-spacing, which provides a scale so that the foliations can be constructed using curves that are nearly linear on a set of some definite size.\\
\item Grim reapers are used in the planar case, but there are no translating solutions on~$\sphere^2$.  Thus we replace them with a suitable solution of a Dirichlet problem.  The analysis is more complicated since the solutions are not explicit.
 \end{enumerate}
 

\tableofcontents{}


\section{A Hausdorff estimate and shrinking circles}

In this section we use the explicit evolution of the shrinking circle to fix a time and scale so that Hausdorff distance between a curve and a great circle will not increase significantly if they are initially close. Lemma~\ref{expand} is used in the sequel to guarantee that a leafable curve does not pass through the endpoints of the solution to the Dirichlet problem constructed in Section~\ref{dirichlet_section}.  We write $B_r(g)$ to denote the $r$-neighbourhood of $g$. 

\begin{lem}~\label{expand} Given $\alpha>0$ there exists $R,T>0$ such that if $g$ is a great circle, $\gamma$ is a smooth closed curve, $0<r<R$ and $0<t<T$, then 
$$
\gamma\subset B_r(g)
$$
implies
$$
\gamma_t\subset B_{(1+\alpha)r}(g).
$$
\end{lem}

We begin by computing the evolution of the shrinking circle since it acts as a barrier for the evolution of $\gamma$.

\begin{lemma}~\label{circles} A shrinking circle with $r_0<\frac{\pi}{2}$ satisfies
$$
r(t)=\arccos\left(\cos(r_0)e^{t}\right).
$$
\end{lemma}

\proof The surface area of a sector with (spherical) radius $r$ is $2\pi(1-\cos(r))$.  Substituting this into the Gauss-Bonnet Formula
$$
\int_{M} K\mathrm{dA}+\int_{\partial M}\kappa_g\mathrm{ds}=2\pi\chi(M)
$$
gives
$$
2\pi(1-\cos(r))+2\pi\sin(r)\kappa_g=2\pi
$$
and hence
$$
\kappa_g=\cot(r).
$$

The Theorem then follows by solving
$$
r'=-\cot(r).
$$\qed

\begin{cor} The extinction time of a circle of radius $r_0<\frac{\pi}{2}$ is 
$\ln(\sec(r_0))$.
\end{cor}


\proof [Proof of Lemma~\ref{expand}] Suppose that $\gamma\subset B_R(g)$, where $R>0$ will be chosen below. Let $\rho_t$ be the radius of the shrinking circle with initial radius $\rho_0=\frac{\pi}{2}-R$.  Then
$$
\gamma_t\subset B_{\frac{\pi}{2}-\rho_t}(g)
$$
since each component of $\partial \mathcal{B}_R(g)$ acts as barrier for the evolution of $\gamma$. 

Defining $R_t=\frac{\pi}{2}-\rho_t$, Lemma~\ref{circles} implies  $R_t=\arcsin(\sin(R)e^t)$.

Now fix $T<\ln(1+\alpha)$ and choose $R>0$ so that
$$
e^T<\frac{\sin((1+\alpha)r)}{\sin(r)}
$$
for all $0<r<R$.  Then for any $0<t<T$ and $0<r<R$
$$
e^t<\frac{\sin((1+\alpha)r)}{\sin(r)}
$$
and hence $R_t<(1+\alpha)r$, as required.\qed


\section{$\mathcal{C}^1$-close to a great circle}

As outlined in the Introduction our argument requires straightening curves so that they are $\mathcal{C}^1$-close to a great circle.  In this section we fix such a notion.

\begin{defn} [Latitudes of $g$] \label{latitudes} Let $g=\partial B_{\frac{\pi}{2}}(x)$ be a great circle.  Then the set of latitudes of $g$ are the curves
$$
\{\partial B_{r}(x)\mid 0<r<\pi\},
$$ 
which we note includes $g$ itself.
\end{defn}

\begin{defn} [$\mathcal{C}^1$ close to $g$] \label{close} Let $g$ be a great circle and let $\gamma$ be a smooth embedded closed curve.   For each $x\in\gamma$, let $u_x$ be the unique latitude of $g$ through $x$ and let 
$\theta_x$ be the angle between $\gamma$ and $u_x$ at $x$.  Then $\gamma$ is $\theta$ $\mathcal{C}^1$-close to $g$ if 
$$
\theta_x\leq\theta
$$
for each $x\in\gamma$.
\end{defn}

\begin{remark} We note that this notion of $\mathcal{C}^1$-close does not imply that the curve in question is $\mathcal{C}^0$-close to $g$ in any sense.  For example, by this definition each latitude of $g$ is 0 $\mathcal{C}^1$-close to $g$.  On the other hand we will only be applying this definition to curve which are already known to be contained in some thin neighbourhood of~$g$.  
\end{remark}

While there are several potential definitions of $\mathcal{C}^1$-close the one property needed in this paper is that any foliation of a neighbourhood $B_r(x)\subset\sphere^2$ by curves which are $\mathcal{C}^1$-close to a great circle $g$ is Bi-Lipschitz equivalent to a similarly straight foliation in $\reals^2$.  

Given $x\in\sphere^2$ and $r>0$ let $B_r=B_r(x)\subset\sphere^2$ and $\widetilde{B}_r=B_r(0)\subset\reals^2$.  Moreover, let $\Pi:\sphere^2\to\reals^2$ be the stereographic projection such that $\Pi(x)=0$.  Then if $\delta_r$ is the appropriate dilation of $\reals^2$
$$
\Phi_r=\delta_r\circ\Pi\mid_{B_r}:B_r\to\widetilde{B}_r
$$
is a conformal diffeomorphism that sends the restriction of great circles through $x$ to line segments through the origin.  If $g$ is a great circle through $x$ then the image of a latitude of $g$ is not a linear segment (unless the latitude is $g$ itself).  However, since the image of the latitudes of~$g$ converge smoothly to $\Phi_r(g)$ we obtain the following:

\begin{lemma} \label{lipschitz} Given $0<\alpha<\beta<\frac{\pi}{2}$  there exists $r=r(\alpha,\beta)$ such that if $g$ is a great circle and $\gamma\subset B_r(g)$ is a smooth closed curve which is $\alpha$ $\mathcal{C}^1$-close to $g$ then $\Phi_r(\gamma\cap B_r)$ is the graph of a $\tan(\beta)$-Lipschitz function over $\Phi_r(g\cap B_r)$.
\end{lemma}

\proof Choose $r>0$ small enough so that if $u$ is a latitude of $g$ then $\Phi_r(u)$ is the graph of a $\tan(\beta-\alpha)$-Lipschitz function.  Then since $\gamma$ is $\alpha$ $\mathcal{C}^1$-close to $g$ and $\Phi_r$ is conformal, $\Phi_r(\gamma)$ makes an angle at most $\beta$ with lines parallel to $\Phi_r(g)$.  This proves the result.\qed 

Lemma~\ref{lipschitz} allows us to obtain the local length estimates in $\sphere^2$ discussed in the Introduction by transferring to a portion of $\reals^2$ and using simple estimates there.


\section{A Straightening Lemma}\label{straight_sec}

In this section we state a Straightening Lemma (Theorem~\ref{straight}) and define the set of curves to which it applies.  The proof is contained in the next section.

The preliminary work is in explaining the choice of various constants.  First, we fix (for the remainder of the paper) a constant $0<\alpha\ll1$.  All other constants depend on $\alpha$ but there is no need to vary it so we suppress this dependence.  One of the roles of $\alpha$ is that the Straightening Lemma straighten curves to be $\alpha$ $\mathcal{C}^1$-close to a great circle in the sense of Definition~\ref{close}.

Next, let $g$ be a great circle and $t>0$.  The goal is to choose a scale $r>0$ and construct a foliation $\mathcal{F}=\{\ell_{\lambda}\}_{\lambda\in[0,1]}$ containing $B_r(g)$ so that the leaves of the evolving foliation, i.e. 
$$
\mathcal{F}_t=\{{(\ell_{\lambda})}_t\}_{\lambda\in[0,1]},
$$
are $\mathcal{C}^1$-close to $g$.  The first step in the proof is a barrier argument.  At $t=0$ the barrier intersects each leaf of the foliation exactly once.   In order to guarantee that this is possible we choose the leaves of $\mathcal{F}$ so as to be controlled on a region $V$ defined below.

Let $0<C<\frac{\pi}{2}$ and $x\in\sphere^2$ and $a:\sphere^2\to\sphere^2$ be the antipodal map.  Define
$$
V=B_C(x)\cup B_C(ax).
$$


Now, let $g$ be a great circle through $x$ and choose $r>0$ satisfying $2r<\alpha C$.  Since $\alpha\ll 1$ this implies that the annulus $B_r(g)$ is thin compared to $V$.  The following definition makes precise the set of curves allowed as leaves of $\mathcal{F}$.  

\begin{defn}\label{leafable} [Leafable] Let $C, x, g, r$ and $V$ be as above.  An embedded smooth closed curve~$\ell$ is {\it leafable} if 
\begin{enumerate}
\item $\ell\subset B_{2r}(g)$, and
\item $\ell\cap V$ is a graph over $V\cap g$ which is $\frac{\alpha}{2}$ $\mathcal{C}^1$-close to $g$.\end{enumerate}
\end{defn}

See Figure~\ref{leafable_fig}.  Note that being leafable implies that $\ell$ is a generator of $\pi_1 B_r(g)$. 

We now state the straightening result that will be used in Section~\ref{foliation_sec}.

\begin{thm}[Straightening Lemma]\label{straight} Given $C>0$, there exists $T>0$ such that for each $0<t<T$ there exists a constant $r=r(t,C)>0$ such that if $\ell\subset B_r(g)$ is leafable then $\ell_t$ is $\alpha$ $\mathcal{C}^1$-close to $g$.  
\end{thm}


We first observe that if $\ell$ is leafable then $\ell_t$ continues to be $\mathcal{C}^1$-close to $g$ on a definite subset of $V$ as long as $t$ and $r$ satisfy Lemma~\ref{expand}.  The argument below is similar to the proof of Theorem~\ref{straight} contained in the next section.  In the proof of Lemma~\ref{maintain} static great circles are used as barriers while the proof of Theorem~\ref{straight} requires specially constructed evolving arcs.

\begin{lem}\label{maintain} Let $C>0$, and let $r>0$ and $t>0$ satisfy Lemma~\ref{expand} and $2r<\alpha C$. Then if~$\ell\subset B_r(g)$ is leafable $$\ell_t\cap B_{\frac{C}{2}}(x)\cup B_{\frac{C}{2}}(ax)$$ consists of two components each of which is $\alpha$ $\mathcal{C}^1$-close to $g$. 
\end{lem}

\proof 

Given $u\in B_{\frac{C}{2}}(x)\cap B_{(1+\alpha)r}(g)$ let $\phi_u$ be the latitude of $g$ containing $u$, and let ${g_1}_u$ and~${g_2}_u$ be the two great circles containing $u$ such that
$$
\angle(\phi_u,{g_i}_u)=\alpha
$$
for $i=1,2$.  Furthermore, define 
$$
\mathcal{G}=\bigcup_{u} \{{g_1}_u,{g_2}_u\},
$$ 
where the union is taken over all $u\in B_{\frac{C}{2}}(x)\cap B_{(1+\alpha)r}(g)$.  

\begin{figure}
\centering
\scalebox{0.6}{\includegraphics{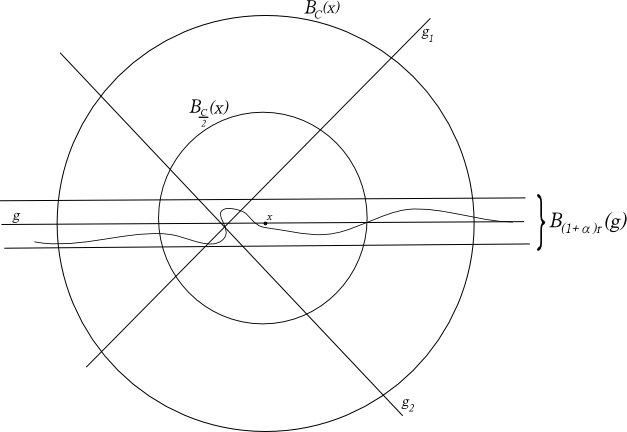}}
\caption{If $\ell_t$ makes an angle greater than $\alpha$ with a latitude of $g$ then it necessarily intersects one of $g_1$ and $g_2$, which are static under the flow, more than once.  Thus Lemma~\ref{maintain} follows from the fact that a leafable curve, which is $\frac{\alpha}{2}$ $\mathcal{C}^1$-close to $g$, intersects each of $g_1$ and $g_2$ exactly once.}  
\label{graph_fig}
\end{figure}

To prove the result it suffices to show that $\ell_t$ intersects each $\rho\in\mathcal{G}$ at most once in each of $B_{\frac{C}{2}}(x)$ and $B_{\frac{C}{2}}(ax)$.  Indeed, suppose that $x\in\ell_t$, $u_x$ is the latitude of $g$ through $x$ and 
$$
\mid\angle({g_i}_{u_x},\ell_t)\mid>\alpha.
$$
Then since $\ell_t$ generates $\pi_1 B_{(1+\alpha)r}(g)$ it follows that $\ell_t$ intersects one of ${g_1}_{u_x}$ or ${g_2}_{u_x}$ a second time.  See Figure~\ref{graph_fig}.

It remains to show that $\ell$ intersects each curve in $\mathcal{G}$ exactly once in $B_C(x)$.  Suppose that there exists $\tilde{g}\in\mathcal{G}$ which intersects $\ell$ more than once in $B_C(x)$.  Now, $2r<\alpha C$ implies that for each $\rho\in\mathcal{G}$
$$
\rho\cap B_{(1+\alpha)r}(g)\subset V,
$$
and hence $\ell_t\cap\rho\subset V$ since Lemma~\ref{expand} implies that $\ell_t\subset B_{(1+\alpha)r}(g)$.    By the mean value argument there exists a latitude of $\tilde{g}$, say $u_{\tilde{g}}$, for which 
$$
\angle(\ell,u_{\tilde{g}})=0
$$
at some point $y$.  Let $u_g$ be the latitude of $g$ through $y$.  Then
$$
\mid\angle(u_g,u_{\tilde{g}})\mid>\frac{\alpha}{2}
$$
since $\tilde{g}$ makes an angle $\alpha$ with some latitude of $g$.

Thus $\mid\angle(\ell,u_g)\mid>\frac{\alpha}{2}$, a contradiction.\qed


\section{Setting up the Dirichlet Problem}\label{dirichlet_section}

The proof of the analogue of Theorem~\ref{straight} in the planar case, which is Lemma 1.7 in~\cite{L}, is a barrier argument using translating solutions to curve shortening flow known as grim reapers.  To create a curve on $\sphere^2$ that serves the same purpose we consider the Dirichlet problem for curve shortening flow, where the initial curve is chosen so that its endpoints lie on either side of a great circle and the interior of the curve follows the geodesic out and back.  We construct the initial curve so that the existence of a well-defined curve shortening flow is guaranteed by the following Theorem of Allen, Layne and Tsukhara:

\begin{thm} [Theorem 1,~\cite{ALT}] \label{existence} Let $\Omega\subset\sphere^2$ be the closed convex domain bounded by two great circles $g_a$ and $g_b$, and let $\Gamma_0:[a,b]\to\sphere^2$ be a smooth embedded curve such that
\begin{enumerate}
\item $\Gamma_0(a,b)\subset \mathrm{int}\Omega,$
\item  $\Gamma_0(a)\in g_a$, 
\item $\Gamma_0(b)\in g_b$, and
\item the curvature and all of its derivatives vanish at $\Gamma_0(a)$ and $\Gamma_0(b)$. 
\end{enumerate}
Then there exists $\Gamma:[a,b]\times[0,\infty)$ with $\Gamma(\cdot,0)=\Gamma_0$ which is a solution to the Dirichlet Problem for curve shortening flow.  Moreover, if $\Gamma_0(a)$ and $\Gamma_0(b)$ are not conjugate then $\Gamma_t$ converges to the unique great circle between them.
\end{thm}

\begin{remark}  Theorem 1 in~\cite{ALT} also applies curves in $\reals^2$ and $\mathbb{H}^2$.  The only change to the statement above is to replace {\it great circle} with {\it geodesic}.
\end{remark}

On $\sphere^2$ convex regions bounded by two great circles are wedges of the form
$$
W_\theta(g,x)=\{R_\psi(g)\mid -\theta\leq\psi\leq\theta\},
$$
where $R_\psi$ is the rotation of $\sphere^2$ by $\psi$ that fixes $x$.  Since $g$ and $x$ are often fixed we write simply $W_\theta$.

\begin{figure}
\centering
\scalebox{0.65}{\includegraphics{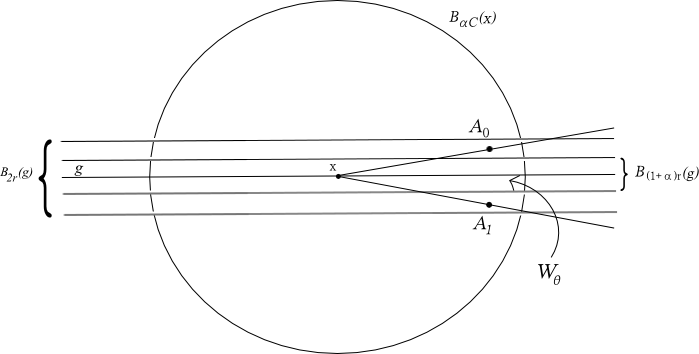}}
\caption{The angle of the wedge $W_\theta$ is chosen so that its boundary exits $B_{\alpha C}(x)$ in $B_{2r}(g)\setminus B_{(1+\alpha)}(g)$.}  
\label{wedge_fig}
\end{figure}

We now define a curve $\Gamma$ whose endpoints lie on a suitably thin wedge~$W_\theta$.  The evolution of $\Gamma$ by curve shortening flow will play the role of the grim reaper.  As in Section~\ref{straight_sec} let $0<C<\frac{\pi}{2}$ and recall that $\alpha\ll1$.  Moreover, we assume that $r$ and $t$ are chosen to satisfy Lemma~\ref{expand} and $2r<\alpha C$.  Since $2r<\alpha C$ there exists $\theta=\theta(r)>0$ such that
\begin{enumerate}
\item $W_\theta\cap B_{\alpha C}(x)\not\subset B_{(1+\alpha)r}(g)$, and
\item $W_\theta\cap \left(B_{\alpha C}(g)\setminus B_{2r}(g)\right)=\varnothing$.
\end{enumerate}

The first condition implies that there exists 
$$
A_0\in \left(B_{\alpha C}(x)\cap B_{2r}(g)\right)\setminus B_{(1+\alpha)r}(g)
$$
which lies on $\partial W_\theta$.  Let $A_1$ be the reflection of $A_0$ across $g$.  These points will be the endpoints of $\Gamma$.  The second condition implies that $\theta(r)\to 0$ as $r\to 0$.   See Figure~\ref{wedge_fig}.


\begin{defn}\label{reaper} For each $r>0$ satisfying $2r<\alpha C$ there exists $\Gamma:[0,1]\to\sphere^2$ with the following properties:

\begin{enumerate}
  \item $\Gamma(0)=A_0$ and $\gamma(1)=A_1$,
  \item $\Gamma(u)\subset B_{2r}(g)\cap W_\theta$ for all $u\in(0,1)$,
  \item $\Gamma\cap B_{(1+\alpha)r}\subset B_{\alpha C}(ax)$, 
  \item $\Gamma$ is a double graph over $g$ (except one point where it intersects $g$),
  \item $\Gamma$ intersects each great circle $R_\psi(g)$ only once, 
   \item $\Gamma\cap B_{\alpha C}(ax)$ makes angle greater than $\frac{\pi}{4}$ with each latitude of $g$, and
   \item the curvature of $\Gamma$ does not change sign, i.e. it is convex.
\end{enumerate}

\end{defn}

\begin{figure}
\centering
\scalebox{0.7}{\includegraphics{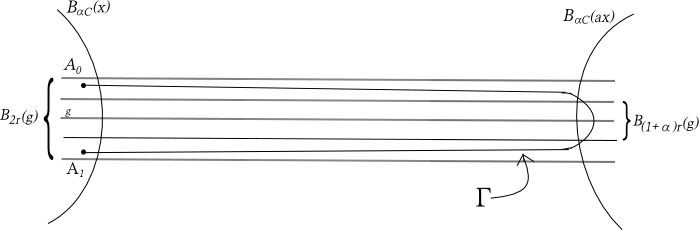}}
\caption{The definition of $\Gamma$.}  
\label{gamma_fig}
\end{figure}

See Figure~\ref{gamma_fig}.

According to Theorem~\ref{existence} we obtain a solution to curve shortening flow $\Gamma_t\subset W_\theta$ that exists for all time and converges to the geodesic between $A_0$ and $A_1$.

Property (6) and Lemma~\ref{maintain} imply that if $\ell$ is leafable and $t>0$ satisfies Lemma~\ref{expand} then $\ell_t$ intersects $\Gamma$ exactly once.  This allows us to 
interpret $\ell_t$ and $\Gamma$ as initial curves and prove the following intersection result:

\begin{lem}\label{shift}  Let $\ell\subset B_r(g)$ be leafable.  Then 
$$
|\ell_t\cap\Gamma_{t^*}|=1
$$
for all $0<t^*\leq t<T$.
\end{lem}

\proof   Let $0<t^*\leq t<T$.  First observe that Lemma~\ref{expand} and Property (3) of Definition~\ref{reaper} imply
$$
\ell_{t-t^*}\cap\Gamma\subset B_{\alpha C}(ax).
$$ 
On the other hand Lemma~\ref{maintain} and Property (6) imply that $\ell_{t-t^*}$ and $\Gamma$ intersect at most once in $B_{\alpha C}(ax)$.  Since $\ell_{t-t^*}$ and $\Gamma$ must intersect by continuity we obtain
$$
|\ell_{t-t^*}\cap\Gamma|=1.
$$
The result follows by evolving each curve for time $t^*$.  The fact that the number of intersections does not increase follows from Lemma~\ref{expand} since it implies that the evolution of $\ell$ does not intersect the endpoints of $\Gamma$, which lie outside $B_{(1+\alpha)r}(g)$.\qed

The next result guarantees that given $t>0$ a scale can be chosen so that the evolution of $\Gamma$ has passed through $B_{2r}(g)\setminus B_\frac{C}{2}(x)$ by time~$t$.

\begin{thm}\label{through} Let $T(r)>0$ be the first time such that $\Gamma_T\subset \overline{B_{\frac{C}{2}}(x)}$.  Then
$$
\lim_{r\to 0} T(r)=0.
$$
\end{thm}

\begin{remark} $T(r)$ exists since the geodesic between $A_0$ and $A_1$ is contained in~$B_{\alpha C}(x)$.
\end{remark}

\proof  Note that $\Gamma_t\subset\overline{B_{\frac{C}{2}}(x)}$ is implied by 
$$
\mathcal{L}(\Gamma_t)<\frac{C}{2}
$$
since the endpoints of $\Gamma$ are contained in $B_{\alpha C}(x)$ and $\alpha\ll 1$. 

Suppose that $\Gamma_t\not\subset B_\frac{C}{2}(x)$.  Then Cauchy-Schwartz implies 
$$
\frac{\partial\mathcal{L}(\Gamma_t)}{\partial t}=-\int_{\Gamma_t}\kappa^2\mathrm{ds}<-\frac{c}{r},
$$
for come $c>0$ and since the length of $\Gamma=\Gamma_0$ is bounded by $3\pi$ this yields
$$
\mathcal{L}(\Gamma_t)<3\pi-\frac{ct}{r}.
$$
Therefore there exists a constant $\tilde{c}>0$ such that $\Gamma_t\subset\overline{B_{\frac{C}{2}}(x)}$ whenever
$$
t>\tilde{c}r
$$
and hence $T(r)<r\tilde{c}$, which proves the result.
\qed

We can now prove Theorem~\ref{straight}.  


\proof [Proof of Theorem~\ref{straight}] Let $t>0$ satisfy Lemmas~\ref{expand} $r=r(t,C)$ be a constant which will chosen later satisfying
\begin{enumerate}
\item $r$ satisfies Lemma~\ref{expand},
\item $2r<\alpha C$, and
\item $T(r)<t$, where $T(r)$ is defined as in Theorem~\ref{through}.
\end{enumerate}

Since $\ell$ is leafable there exists $x\in\mathbb{S}^2$ such that $\ell$ is $\frac{\alpha}{2}$ $\mathcal{C}^1$-close to $g$ on $B_C(x)$, and Lemma~\ref{maintain} implies the conclusion of the Theorem is true on $U:=B_{\frac{C}{2}(x)}\cup B_{\frac{C}{2}(ax)}$  


Suppose that the endpoints of $\Gamma$  lie on the wedge $W_{\theta_0}$, where $\theta_0$ will be chosen below.  By construction $\Gamma$, and hence $\Gamma_t$, intersect each great circle $R_\theta(g)$ exactly once for all $t>0$ and $|\theta|<\theta_0$.  Hence the function $A:[-\theta_0,\theta_0]\times[0,\infty)\to (0,\pi)$ defined by
$$
A(\theta,t)=\angle(\Gamma_t,R_\theta(g))
$$
is well-defined.

By further choosing $r$, and hence $\theta_0$, small enough it is possible to arrange that 
$$
A(\theta_0,t)<\alpha/8
$$ 
whenever $t<T_r$, i.e. $\Gamma_t\not\subset B_{\frac{C}{2}}(x)$.  This can be done since the convexity of $\Gamma$ implies that $A(\cdot,t)$ is monotonic.  Indeed $\Gamma_t\not\subset B_{\frac{C}{2}}(x)$ implies that $\Gamma_t$ follows $g$ for some definite distance while the angle between $R_{\theta_0}(g)$ and $g$ tends to zero as $\theta_0\to 0$. 
  
Now, since $A$ is continuous there exists $\epsilon>0$ such that $|A(\theta,t)|<\alpha/4$ whenever $t<T_r$ and $\theta\in[\theta_0-\epsilon,\theta_0]$.  Thus defining $\Gamma_R=R_{\theta_0-\epsilon}(\Gamma)$, i.e. by rotating $\Gamma$, we have
$$
\angle((\Gamma_R)_t,g)<\alpha/4
$$
for all $t<T(r)$.  Again, by continuity there exists $\theta_1>0$ such that 
$$
\angle((\Gamma_R)_t,R_\psi(g))<\alpha/2
$$
for all $t<T(r)$ and $|\psi|<\theta_1$.  


As in the proof of the planar case we use not only $\Gamma_t$ and its rotations but also mirror images of which pass through $B_r(g)$ in the opposite direction.  Let $\tilde{g}$ be the great circle which is constant distance from $x$ (and~$ax$), and define $\tilde{\Gamma}_R$ to be the reflection of $\Gamma_R$ across~$\tilde{g}$.

Now suppose that $z\in\ell_{t}\cap B_{(1+\alpha)r}(g)\setminus U$.  Since $T(r)<t$ there exists $t_1,t_2<t$ such that $z\in({\Gamma_R})_{t_1}\cap{(\widetilde{\Gamma}_R)}_{t_2}$.  To simplify notation we write $\Gamma_1=({\Gamma_R})_{t_1}$ and $\widetilde{\Gamma}_2={(\widetilde{\Gamma}_R)}_{t_2}$.  Then Lemma~\ref{shift} implies
$$
z=\Gamma_1\cap\ell_t=\widetilde{\Gamma}_2\cap\ell_t.
$$
Let $Z$ be the the convex hull of $\Gamma_1$ and $\widetilde{\Gamma}_2$.  Since $\ell_t$ intersects each of $\Gamma_1$ and $\widetilde{\Gamma}_2$ once it follows that $\ell_t$ cannot leave $Z$ outside $U$ since then returning to $Z$ would cause a second intersection.  Thus, if $g_z$ is the latitude of $g$ containing $z$ we have
$$
-\alpha\leq\angle(\widetilde{\Gamma}_2,g_z)\leq\angle(\ell_t,g_z)\leq\angle(\Gamma_1,g_z)\leq\alpha
$$
as required.\qed


\section{$(C,\theta)$-spacings}\label{space_sec}

In order to prove Theorem~\ref{length} it is necessary to produce, for each $x\in\sphere^2$, two foliations through $x$ (with leafable leaves) that evolve to be nearly perpendicular.  In this section we show that for each Jordan curve there is a scale so that this can be done.  Here, scale refers to the constant $C>0$ in the definition of leafable.

In what follows we denote the antipodal map by $a$ and write $\overline{xy}$ for the unique great circle containing non-antipodal points $x$ and $y$.  

\begin{defn}\label{spacing} [$(C,\theta)$-spacing] Given a Jordan curve  $\gamma$ on $\sphere^2$ and constants $C>0$ and $0<\theta<\frac{\pi}{2}$ we say that a set of points $\{y_1,y_2,\ldots, y_n\}$ is a {\it $(C,\theta)$-spacing} for $\gamma$ if
\begin{enumerate}
\item $\overline{B_C(y_i)\cup B_C(ay_i)}\cap\gamma=\varnothing$ for each $i=1,\ldots,n$
\item For each $x\in\sphere^2$ there exists $i_1,i_2$ such that the angle between the great circles 
$\overline{xy_{i_1}}$ and $\overline{xy_{i_2}}$ is greater than $\frac{\pi}{2}-\theta$.
\end{enumerate}
\end{defn}

\begin{lemma} \label{ex_space} Let $\gamma$ be a Jordan curve on $\sphere^2$. For each $0<\theta<\frac{\pi}{2}$ there exists $C>0$ such that a $(C,\theta)$-spacing exists for $\gamma$.
\end{lemma}

\proof Let $x\in\sphere^2$ and let $u_1$, $u_2$ be points satisfying

\begin{enumerate}
\item $d(x,u_i)=\frac{\pi}{2}$ for each $i=1,2$, and
\item $d(u_1,u_2)=\frac{\pi}{2}$.
\end{enumerate}

Then $\overline{xu_1}$ and $\overline{xu_2}$ are perpendicular and by continuity there exists $\epsilon_x=\epsilon_x(\theta)>0$ such that if $v_i\in B_{\epsilon_x}(u_i)$ for $i=1,2$ and $w\in B_{\epsilon_x}(x)$ then the angle between $\overline{wu_1}$ and $\overline{wu_2}$ is at least $\frac{\pi}{2}-\theta$. 

For $i=1,2$ choose $y^x_i\in B_{\epsilon_x}(u_i)$ so that $y_i, ay_i\notin\gamma$.  The existence of such points is a consequence of the Jordan-Schoenflies Separation Theorem.  Including $y^x_1$ and $y^x_2$ in a potential spacing for $\gamma$ guarantees that Condition (2) of Definition~\ref{spacing} is satisfied for all points in $B_\epsilon(x)$.

Repeating this process at each $x$ we obtain a cover of $\sphere^2$ by set of the form $B_{\epsilon_x}(x)$.  Let $\{x_i\}$ be a finite set such that $\{B_{{\epsilon_x}_i}(x_i)\}$ is a finite subcover.  Then the finite set $\bigcup\{y^{x_i}_1,y^{x_i}_2\}$ satisfies Condition (2).

Finally, since no point of $\bigcup\{y^{x_i}_1,y^{x_i}_2\}$  lies on $\gamma$ it is a $(C,\theta)$-spacing for some $C>0$.\qed      

When $\gamma$ is smooth Theorem~\ref{length} [provides length estimates for the evolution of $\gamma$ on a finite time interval $(0,T)$ and $T$ depends only on the values of $C$ and $\theta$ for which a $(C,\theta)$-spacing exists.  Together with the observation below this implies that if $\gamma_n$ is a uniformly converging sequence then for sufficiently large $n$ the length estimates are valid on a uniform interval. 

\begin{lemma}\label{spacing} Let $\gamma_n$ be a sequence of Jordan curves that converge uniformly to a Jordan curve $\gamma$.  Then any $(C,\theta)$-spacing for $\gamma$ is a $(C,\theta)$-spacing for $\gamma_n$ for sufficiently large $n$.
\end{lemma} 

\proof Condition (2) of Definition~\ref{spacing} does not depend on the curve and Condition~(1) will be satisfied for sufficiently large $n$ since
$$
\bigcup_{i}\overline{B_C(y_i)\cup B_C(ay_i)}
$$
is compact. \qed


\section{$r$-multiplicity}\label{mult_sec}

In this section we define the $r$-multiplicity.  It is a straightforward generalization of the idea of the same name that appeared in \cite{L}  and more detail can be found there.  For example, \cite{L} contains a compactness result for sets of curves that satisfy an $r$-multiplicity bound at all scales.  In \cite{L} the case of immersed curves was considered but here we restrict to the embedded case since that is all that is needed in the proof of Theorem~\ref{JordanCurve}.  This greatly simplifies the exposition.  In what follows a Jordan curve is a continuous embedding of $\sphere^1$ into $\sphere^2$ but we often identify the curve with its image.    

\begin{defn} [$r$-multiplicity ] Let $g$ be a great circle, $0<r<\frac{\pi}{2}$ and $\gamma$ be a Jordan curve in $\mathbb{S}^2$.  Then the $r$-multiplicity of $\gamma$ at $g$, denoted by $M_{r,g}(\gamma)$, is defined as the number of components of $\gamma\cap B_{2r}(g)$ which intersect $\overline{B_r(g)}$ non-trivially. 

In addition we define {\it the $r$-multiplicity of $\gamma$} by
 $$M_r(\gamma)=\sup_g\{M_{r,g}(\gamma)\}.
 $$  
\end{defn}

The $r$-multiplicity acts as a coarse intersection number.  The coarseness is important; in contrast to the fact that a Jordan curve may intersect a straight line infinitely many times we have the following:

\begin{lemma}\label{mult_converge} Let $\gamma$ be a Jordan curve.  Then for $M_r(\gamma)<\infty$ for any $r>0$.
\end{lemma}

\proof The result follows from the uniform continuity of $\gamma:\mathbb{S}^1\to\mathbb{S}^2$.\qed

Moreover, since the $r$-multiplicity implicitly involves counting the number of intersections with the static lines the fact  that the number of intersections between evolving curves does not increase~\cite{An88} can be used to show the $r$-multiplicity is monotonic.

\begin{lemma}\label{multmonotonic}  $M_{r,g}(\gamma_t)$ is non-increasing in $t$. 
\end{lemma}

It is also straightforward to verify that the $r$-multiplicity behaves well under uniform convergence:

\begin{lemma} \label{converge} Let $\gamma_n$ be a sequence of Jordan curves that converges uniformly to a Jordan curve $\gamma$.  Then for each great circle $g$ and $r>0$
$$
\limsup_{n\to\infty}M_{r,g}(\gamma_n)\leq M_{r,g}(\gamma).
$$ 
\end{lemma}

In Lemma~\ref{converge} the inequality comes from the case where the $M_{r,g}(\gamma)$ counts a component of $\gamma\cap B_{2r}(g)$ that intersects $\overline{B}_r(g)$ but not $B_r(g)$.  In this case the corresponding arcs in $\gamma_n$ need never intersect $\overline{B}_r(g)$.

We note that Theorem~\ref{approx} is an immediate consequence of Lemmas~\ref{ex_space} and~\ref{mult_converge} and Theorem~\ref{length}.


\section{Foliations}\label{foliation_sec}

The proof of Theorem~\ref{length} proceeds by constructing a foliation $\F$ of the annulus $B_r(g)$ with the property that each leaf of the foliation intersects a given curve $\gamma$ at most $2M_{r,g}(\gamma)$ times.  Recall that $V=B_C(y)\cup B_C(ay)$.  If $y$ is a point in a $(C,\theta)$-spacing for $\gamma$ then
$\gamma\cap V=\emptyset$ implying that the definition of the leaves of $\mathcal{F}$ on $V$ does not affect the number of intersections with $\gamma$.  In particular this guarantees that $\mathcal{F}$ can be constructed so that each leaf is leafable.

We now construct the initial foliation.  The construction is essentially the same as in~\cite{L} except here the exposition is considerably simpler since we consider only embedded curves.  

\begin{thm} [The initial foliiation]\label{foliate1} Let $C, \theta>0$ and let $r>0$ satisfy $2r<\alpha C$ and Lemma~\ref{expand}.  Let $\gamma$ be a smooth closed curve for which there exists a $(C,\theta)$-spacing, and let $y$ be a point in such a spacing.  Then for each  great circle $g$ containing $y$ there exists a 1-parameter family of smooth curves
$\mathcal{F}=\{\ell_x\}_{x\in[0,1]}$ such that:
\newline\indent (1) $\mathcal{F}$ foliates a region containing the annulus $B_r(g)$.
\newline\indent (2) $\ell_x$ is leafable for each $x\in[0,1]$, and
\newline\indent (3) $|\gamma\cap\ell_x|\leq 2M_{r,g}(\gamma)$ for
each $x\in [0,1]$.
\end{thm}

\proof The first step is to define $\ell_0$, which lies in a component of $B_{2r}(g)\setminus B_r(g)$.  Since the definition of $(C,\theta)$-spacing implies that $\gamma\cap B_C(y)=\emptyset$ we define $\ell_0$ to coincide with a latitude of $g$ in $B_C(y)$.   

On $B_{2r}(g)\setminus B_C(y)$ define $\ell_0$ so that it intersects each component $\tilde\gamma$ of  $\gamma\cap B_{2r}(g)$ transversely according to the following scheme.  Note this is the minimum number of intersections necessary if $\ell_0$ is to remain in $ B_{2r}(g)\setminus B_C(y)$.
\begin{enumerate}
    \item $|\tilde\gamma\cap\ell_0|=0$ if $\tilde\gamma$ does not count towards $M_{r,g}(\gamma)$,
    \item $|\tilde\gamma\cap\ell_0|=1$ if $\tilde\gamma$ contributes to $M_{r,g}(\gamma)$ and the endpoints of $\gamma$ lie on the same component of $\partial B_{2r}(g)$, and
    \item $|\tilde\gamma\cap\ell_0|=2$ if $\tilde\gamma$ contributes to $M_{r,g}(\gamma)$ and the endpoints of $\gamma$ lie on the distinct components of $\partial B_{2r}(g)$.
\end{enumerate}
Similarly, one defines $\ell_1$ in the second component of $B_{2r}(g)\setminus B_r(g)$.

Now, let $L$ be the closed annulus between $\ell_0$ and $\ell_1$.   Then $\gamma\cap L$ contains $M=M_{r,g}(\gamma)$ arcs which intersect $\ell_0$ and $\ell_1$ transversely.  Let $S=\sphere^1\times[0,1]$ and let $\{\mu\}_{i=1}^M$ be a collection of arcs in $S$ either of which is the restriction of a vertical line or the graph of a parabola arranged so the combinatorial structures of $(L,\gamma\cap L)$ and $(S,\{\mu_i\})$ are equivalent.  Then there exists a diffeomorphism $\Phi:S\to L$ which sends each arc in $\gamma\cap L$ to an arc in  $\{\mu_i\}$.  The foliation
$$
\mathcal{F}=\{\Phi(\{y=\zeta\})_{\zeta\in [0,1]}
$$
then satisfies the requirements of the Theorem.\qed


With $\F$ now defined let $$\F_t=\{{(\ell_x)}_t\}_{x\in[0,1]}$$ be the result of evolving each leaf of $\F$ by curve shortening flow for time $t$.  Since each leaf of $\F$ is leafable we simultaneously apply the Straightening Lemma to each leaf:

\begin{thm} \label{foliate2} Given $t>0$ satisfying Lemma~\ref{expand} and $C>0$ there exists $r=r(t,C)$ such that if $\F$ is constructed as in Lemma~\ref{foliate1} at scale $r$ then
\newline\indent (1) $\mathcal{F}_t$ foliates a region containing the annulus $B_r(g)$,
\newline\indent (2) $(\ell_x)_t$ is $\alpha$ $\mathcal{C}^1$-close to $g$ for each $x\in[0,1]$, and
\newline\indent (3) $|\gamma_t\cap(\ell_x)_t|\leq 2M_{r,g}(\gamma)$ for
each $x\in [0,1]$.
\end{thm}

\proof For (1) note that the evolution of $\partial B_r(g)$ acts as a barrier for $(\ell_0)_t$.  Property (2) follows immediately from Theorem~\ref{straight} and (3) follows from the fact that the number of intersections do not increase under curve shortening flow. \qed

The fact that the foliation $\mathcal{F}_t$ consists of curves $\mathcal{C}^1$-close to a common great circle $g$ does not immediately imply that there is a uniform bound on the biLipschitz constant needed to map $\mathcal{F}_t$ to a set of latitudes of $g$.  Nevertheless, as in~\cite{L} we use the fact that the leaves have been evolving by curve shortening flow to establish their uniform separation.

 The main observation is that the derivative of the holonomy map is proportional to the solution of the linearization of curve shortening flow.  The following is then a consequence of the Harnack inequality.  

\begin{thm}\label{harnack} There exists a constant $d>0$ depending only on $r$ and $t$ such that the foliation $\mathcal{F}_t$ is $d$-bi-Lipschitz equivalent to the foliation of $B_r(g)$ by latitudes of $g$.
\end{thm}
 
The argument is exactly the same as in~\cite{L} to which refer the reader for the details.    But here we do explain what allows us to apply the Harnack inequality in this new setting.  In~\cite{Hsu} Hsu showed that the evolution of a graph over a great circle can be computed by projecting to a Euclidean equation. That is, let
$$
\hat\sphere^2=\{(x,y,z)\mid x^2+y^2+z^2=1, z\neq\pm 1)\},
$$
$$
\mathcal{C}=\{(x,y,z))\mid x^2+y^2=1\}
$$
and $\Pi:\hat{\sphere}^2\to\mathcal{C}$ be the natural radial projection.  Then any $2\pi$-periodic function determines a curve on $\mathcal{C}$ in a natural way and we have the following.  See~\cite{Hsu} for the computations:

\begin{thm} [\cite{Hsu}] If $u_t$ satisfies
\begin{equation}\label{csfsphere} 
u_t=\frac{(1+u^2)^2}{1+u^2+u_x^2}(u_{xx}+u),
\end{equation}
then $\Pi^{-1}(u_t)$ satisfies curves shortening flow.
\end{thm}

In~\cite{Hsu} the linearization of (\ref{csfsphere}) about a solution $u$ was shown to be
$$
v_t=(1+a(u,u_x))v_{xx}+b(u,u_x)v_x+(1+c(u,u_x))v,
$$
where $a$, $b$ and $c$ are small when $u$ and $u_x$ are, and hence solutions to such an equation satisfy he Harnack inequality.


\section{Proof of the length estimate}


\proof [Proof of Theorem~\ref{length}] Recall that $0<\alpha\ll 1$ is a fixed constant.  Let $T<\ln(1+\alpha)$ so that Lemma~\ref{expand} holds.  Let $\gamma$ be a Jordan curve with a $(C,\theta)$-spacing for some $\theta\ll\alpha$.  For each $0<t<T$ we choose the scale $r=r(t,\theta, C)$ such that
\begin{enumerate}
\item the pair $(r,t)$ satisfies Lemma~\ref{expand},
\item $r<r(\alpha,2\alpha)$ in Lemma~\ref{lipschitz},
\item $2r<\alpha C$, 
\item $2T(r)<t$, where $T(r)$ is the function in Theorem~\ref{reaper}, and  
\item $r$ satisfies the conclusion of Theorem~\ref{straight}.
\end{enumerate}

Let $x\in\sphere^2$ and let $y_1$ and $y_2$ be two points in the $(C,\theta)$-spacing such that
$$
\angle(\overline{xy_1},\overline{xy_2})>\frac{\pi}{2}-\alpha.
$$
Applying Theorem~\ref{foliate1} twice we obtain a pair of foliations $\F_1$ and $\F_2$ associated to $\overline{xy_1}$ and $\overline{xy_2}$ respectively.  Conditions (4) and (5) above imply that each $(\F_i)_{t/2}$ satisfies the conclusion of Theorem~\ref{foliate2}.

After allowing the graphical foliations to evolve further on $[t/2,t]$ Lemma~\ref{lipschitz} and Theorem~\ref{harnack} imply that there exists a d-bi-Lipschitz homeomorphism 
$$\Phi_r:B_r(x)\subset\sphere^2\to B_r(0)\subset\reals^2$$ 
which sends the pair of foliations ${\F_i}_t\cap B_r(x)$ to the standard grid in $B_r(0)\subset\reals^2$. . Moreover, since the intersection numbers are preserved by the homeomorphism we obtain a perhaps disjoint curve $\Phi(\gamma_t\cap B_r(x))\subset B_r(0)\subset\reals^2$ that intersects each horizontal and vertical line at most $2M_{r,g}(\gamma)$ times.  A relatively simple calculus exercise yields
$$
\L(\Phi(\gamma_t\cap B_r(x)))<8rM_r(\gamma),
$$
and hence
$$
\L(\gamma_t\cap B_r(x))<8rdM_r(\gamma).
$$
Since $\sphere^2$ can be covered by $\widetilde{C} r^{-2}$ balls of radius $r$ we obtain the estimate
$$
\L(\gamma_t)<\frac{8\widetilde{C}d}{r}M_r(\gamma).
$$
This completes the proof.\qed


\section{Smoothness of the level-set flow}

\subsection{Level-set flow}\label{lsfdef_section}

The level-set flow of a compact set is defined by the property of being the largest evolution that satisfies the avoidance principle.  When the initial data is smooth the level-set flow and curve shortening flow agree until the latter ceases to exist.  The geometric version used here, framed in terms of weak-set flows, was first developed by Ilamanen~\cite{I94}, while the original analytic viewpoint was developed independently in~\cite{CGG91} and~\cite{ES91}.  See also~\cite{HO,I94,W00,W03,LamLauer}.  Here we specialize to $\sphere^2$.

\begin{defn}\label{weaksetflow}\emph{(Weak-set flow, Level-set flow)}
Let $K \subset \sphere^2$ be compact, and let ${\{K_t\}}_{t \geq 0}$ be a 1-parameter family of compact sets with $K_0 = K$, such that the space-time track $\cup (K_t \times \{t\}) \subset \sphere^2\times\reals$ is closed. Then  ${\{K_t\}}_{t \geq 0}$ is a \emph{weak-set flow} for $K$ if for every smooth curve shortening flow $\gamma_t$ defined on $[a,b] \subset [0, \infty]$ we have 
$$
K_a \cap \gamma_a = \emptyset \Longrightarrow K_t \cap \gamma_t = \emptyset
$$
for each $t \in [a,b]$.

The \emph{level-set flow} of a compact set $K \subset\sphere^2$ is the maximal weak-set flow. That is, a weak set flow $K_t$ such that if $\widehat{K}_t$ is any other weak set flow, then $\widehat{K}_t \subset K_t$ for all $t \geq 0$.
\end{defn}

In our case there is an explicit description:  Let $\gamma$ be a Jordan curve and let $\gamma_t$ be its level-set flow.  Let~$\Omega_n$ be an exhaustion of one component of $\sphere^2\setminus\gamma$ by smooth disks, and define $\alpha^n=\partial\Omega_n$.  By repeating this procedure in the other component of $\sphere^2\setminus\gamma$ we obtain a second sequence $\beta^n$.   Now, let $A^n$ be the sequence of nested annuli between ${\alpha^n}$ and ${\beta^n}$ and let $A^n_t$ be the annulus between $\alpha^n_t$ and $\beta^n_t$, the time $t$ evolutions of $\alpha^n$ and $\beta^n$ by curve shortening flow.  Then $\gamma\subset A^n$ and the avoidance principle implies that for any $t>0$
\begin{equation}\nonumber
\gamma_t\subset\bigcap_n A^n_t.  
\end{equation}
And since it is easy to verify that the right-hand side is in fact a weak-set flow we obtain
\begin{equation}\label{lsf}
 \gamma_t=\bigcap_n A^n_t.   
\end{equation}
As a consequence of this, although it also easy to verify directly, we see that the set $\cap A^n_t$ does not depend on the original choice of approximating curves.


\subsection{Proof that the level-set flow is smooth}  With Theorem~\ref{length} now established the smoothness of the level-set flow is proved exactly as in~\cite{L} to which we refer the reader for details. Here we give an outline.

\begin{thm}\label{smooth} (Smoothness) Let $\gamma$ be a Jordan curve on $\sphere^2$ with level-set flow $\gamma_t$.  Then for $t>0$ sufficiently small either
\begin{enumerate}
    \item $\gamma_t$ is a smooth closed curve, or
    \item $\gamma_t$ is an annulus with smooth boundary.
\end{enumerate}
\end{thm}

\proof [Outline of the proof of Theorem~\ref{smooth}] Let $\Omega$ be a domain in $\sphere^2$ such that $\gamma=\partial\Omega$ is a Jordan curve.  Let $\Omega_n\subset\Omega$ be a sequence of smooth disks which exhaust $\Omega$.  This implies that the smooth curves $\alpha^n=\partial\Omega_n$ Hausdorff converges to $\gamma$. 

Theorem~\ref{length} and Lemma~\ref{converge} imply that for each $t>0$ there exists $C_0(t)$ such that
\begin{equation}\label{lengthbound}
\mathcal{L}(\alpha^n_t)<C_0(t). 
\end{equation}

If we denote by $(\Omega_n)_t$ the region bounded by ${\gamma_n}_t$ and 
$$
\Omega_t=\cup(\Omega_n)_t
$$
then (\ref{lsf}) implies that $\partial\Omega_t\subset\partial\gamma_t$ and (\ref{lengthbound}) implies that $\mathcal{H}^1(\partial\Omega_t)<\infty$, where $\mathcal{H}^1$ is the one-dimensional Hausdorff measure.

Recall that for a smooth curve $\alpha_t$ evolving by curve shortening flow
$$
\frac{\mathrm{d}\mathcal{L}(\alpha_t)}{\mathrm{dt}}=-\int_{\alpha_t}\kappa^2\mathrm{ds}.
$$
By comparing $\mathcal{L}(\alpha^n_t)$ at two positive times $0<t_1<t$ we obtain a second constant $C_1(t)>0$ such that

\begin{equation}\label{curvaturebound}
\int_{\alpha^n_t}\kappa^2\mathrm{ds}<C_1(t)
\end{equation}

Together~\ref{lengthbound} and~\ref{curvaturebound} imply that $\partial{\Omega}_t$ is a $\mathcal{C}^1$-curve and since $\partial{\Omega}_t$ is the boundary of the level set flow of $\gamma$ this implies that in fact $\partial{\Omega}_t$ is smooth.  \qed


\section{Uniqueness}\label{unique_section}

From the previous section we know that for small positive times the level-set flow of a Jordan curve is either a smooth closed curve or the region between two disjoint smooth curves.  In this section we show that only the former case occurs when the initial data has measure zero, establishing the uniqueness portion of Theorem 1.1.   We also explore the possible outcomes when the initial data has positive area.  

In Section 2 the Gauss-Bonnet Theorem was used to compute the explicit evolution of a shrinking circle.  It was also used by Gage~\cite{Gage} to show that for a smooth curve on $\sphere^2$ the property of bisecting the area is preserved.  Later it was shown that in fact any such curve converges to a unique great circle.    

Here we again use the Gauss-Bonnet Theorem, this time to compute the change in area of an evolving annulus.  We then apply this computation to the annuli defining the level-set flow as per the discussion preceding (\ref{lsf}).  Besides the Gauss-Bonnet Theorem we also use the following fact which appears as Lemma 1.3 in~\cite{Gage}. 

\begin{lem} \label{Gage}  Let $\gamma_t$ be a closed curve evolving by curve shortening flow on $\sphere^2$.  Then
$$
\frac{d}{dt}\int_{\gamma_t}\kappa_g\mathrm{ds}=\int_{\gamma_t}\kappa_g\mathrm{ds}.
$$
\end{lem}

We use $\mu$ to denote the 2-dimensional Lebesgue measure on $\sphere^2$.

\begin{thm}\label{uniqueness}Let $\gamma$ be a Jordan curve with $\mu(\gamma)=0$.  Then $\mu(\gamma_t)=0$ for each $t>0$.  
\end{thm}

\proof As in Section~\ref{lsfdef_section} let $\alpha^n$ and $\beta^n$ be sequences or smooth approximations that define the level-set flow, and let $A^n_t$ be the annulus bounded by their evolutions by curve shortening flow.  Fix the unit normal on $\alpha^n$ and $\beta^n$ which is outward and inward pointing respectively with respect to $A^n$.  Then the Gauss-Bonnet Theorem implies
\begin{eqnarray}
\nonumber \mu(A^n_t)&=&\int_{\alpha^n_t}\kappa_g\mathrm{ds}-\int_{\beta^n_t}\kappa_g\mathrm{ds},
\end{eqnarray}
and thus
\begin{eqnarray}
\nonumber\frac{d}{dt}\mu(A^n_t)&=&\frac{d}{dt}\left[\int_{\alpha^n_t}\kappa_g\mathrm{ds}-\int_{\beta^n_t}\kappa_g\mathrm{ds}\right]\\\nonumber
    &=&\int_{\alpha^n_t}\kappa_g\mathrm{ds}-\int_{\beta^n_t}\kappa_g\mathrm{ds}\\\nonumber
    &=&\left(-\mu(\alpha^n_t)+2\pi\right)-\left(-\mu(\beta^n_t)+2\pi\right)\\\nonumber
    &=&\mu(A^n_t),
\end{eqnarray}
where Lemma~\ref{Gage} is used at the second equality and Gauss-Bonnet is used at the third.  This equation is valid for all $t>0$ such that $\alpha^n_t$ and $\beta^n_t$ exist.  

It follows from~(\ref{lsf}) that 
$$
\frac{d}{dt}\mu(\gamma_t)=\mu(\gamma_t)
$$
which proves the result since $\mu(\gamma_0)=0$.\qed


The next result asserts that the long-term behaviour of a measure zero Jordan curve satisfies the exact same dichotomy as smooth curves.  Let $\mathrm{Area}(\gamma)\in(0,2\pi]$ denote the area contained in the interior of $\gamma$. 

\begin{thm} If $\mu(\gamma)=0$ and $\gamma$ bisects the area of $\sphere^2$ then $\gamma_t$ is non-empty for all time.  Otherwise, the evolution becomes extinct in finite time.
\end{thm}

\proof If $\gamma$ bisects the area then 
$$
\lim_{n\to\infty}\mathrm{Area}(\alpha^n)=2\pi,
$$
and hence by the Gauss-Bonnet Theorem
$$
\lim_{n\to\infty}\int_{{\alpha^n}}\kappa_g\mathrm{ds}=0.
$$
On the other hand Lemma~\ref{Gage} implies that
$$
\mathrm{Area}(\alpha^n_t)+e^t\int_{{\alpha_n}}\kappa_g\mathrm{ds}=2\pi
$$
for sufficiently small $t$ and hence
$$
\lim_{n\to\infty}\mathrm{Area}(\alpha^n_t)=2\pi.
$$
Since this argument applies equally well to an approximating sequence in either component of $\sphere^2\setminus\gamma$ it follows that~$\gamma_t$ bisects the area for small positive times.  But since $\gamma_t$ is smooth for $t>0$ the result of Gage~\cite{Gage} that smooth curves continue to bisect the area completes the proof.

The proof of the second statement follows by comparing with a curve $\zeta$ with $\mathrm{Area}(\zeta)<2\pi$ and having the property that $\gamma$ lies in the component of its complement with least area.\qed

For positive area curves the computation carried out in the proof of Theorem~\ref{uniqueness} can be used to show that the area does not decrease, and hence that the evolution is a smooth annulus for small positive times. We summarize the long-term behaviour below.  The proofs are left to the reader.

\begin{thm}\label{positive} Let $\gamma$ be a positive measure Jordan curve, $U_1$ and $U_2$ be the two components of $\sphere^2\setminus\gamma$, and
$$
A=\max\{\mathrm{Area}(U_1),\mathrm{Area}(U_2)\}.
$$
Then $\gamma_t$ is initially an evolving annulus with smooth boundary components.
\begin{enumerate}
    \item If $A>2\pi$ then $\gamma_t$ becomes extinct in finite time.
    \item If $A=2\pi$ then one boundary component converges to a great circle and $\gamma_t$ converges to a hemisphere.
    \item If $A<2\pi$ then $\gamma_t=\sphere^2$ for sufficiently large $t$.
    
\end{enumerate}
\end{thm}


\section{Convergence}


The only backward convergence that is guaranteed by the definition of level-set flow is Hausdorff convergence.  In this section we show that as $t\to 0^+$ the level set flow of $\gamma$ converges to the initial data in the $\mathcal{C}^0$-metric on the space of unparametrized curves.  The proof is exactly the same as in the planar case.  The idea is that the existence of `converging overlaps' would contradict Lemma~\ref{multmonotonic}, the fact that $r$-multiplicity is non-increasing.

\begin{thm} $\lim_{t\to 0^+} \gamma_t=\gamma$ in the space of (unparametrized) continuous curves.
\end{thm} 

\proof Suppose otherwise.  Then there exists a sequence $t_i\to 0^+$ and distinct points $x_i,y_i, z_i\in \gamma_{t_i}$ such that 
\begin{enumerate}
\item $A:=\lim_{t\to 0^+} x_i=\lim_{t\to 0^+} z_i\neq \lim_{t\to 0^+} y_i=:B$, and
\item  the two sequences of arcs, one between $x_i$ and $y_i$ and the other between $y_i$ and $z_i$ Hausdorff converge to the same arc in $\gamma$ between $A$ and $B$. 
\end{enumerate}

Let $g$ be a great circle that separates $A$ and $B$ and let $r>0$ be chosen so that $B_r(g)$ does not contain either $A$ or $B$.  We claim that it is possible to choose $g$ and $r$ such that each component defining $M_{r,g}(\gamma)$ enters $B_r(g)$ and not just $\overline{B_r}(g)$.  To see this fix $r>0$ small and use the uniform continuity of $\gamma$ to bound the number of $g$ for which the claim fails to hold.

The above claim implies that if $\alpha_n$ is any sequence Hausdorff converging to $\gamma$ then
$$
M_{r,g}(\alpha_n)\geq M_{r,g}(\gamma)
$$
for $n$ sufficiently large.  In our case property (2) above implies that the inequality is strict, that is
$$
M_{r,g}(\gamma_t)>M_{r,g}(\gamma)
$$
for sufficiently small $t>0$.

But then since $\gamma_t$ Hausdorff converges to $\gamma$ it follows that when $t$ is small $\gamma_t$ also has the property that each component defining $M_{r,g}(\gamma_t)$ intersects $B_r(g)$ non-trivially.  Hence we have equality in Theorem~\ref{converge} and so
 $$
 M_{r,g}((\gamma_n)_{t_i})>M_{r,g}(\gamma)
 $$ 
 for $n$ sufficiently large and $t_i$ sufficiently small.   Finally by Theorem~\ref{multmonotonic} the $r$-multiplicity is non-increasing and so 
 $$
 M_{r,g}(\gamma_n)>M_{r,g}(\gamma),
 $$
 contradicting Theorem~\ref{mult_converge}.  This completes the proof.\qed



Fachbereich Mathematik und Informatik, Freie Universit$\ddot{a}$t, Berlin, Germany, 14195\\
\phantom{tt} Email address: lauer@zedat.fu-berlin.de

\end{document}